\documentclass[12pt]{article}
\usepackage{indentfirst}
\usepackage{amsmath}
\usepackage{amssymb}

\newtheorem{Th}{Theorem}
\newtheorem{Lem}{Lemma}

\title{The distribution of second degrees in the Bollob\'as--Riordan random graph model
\footnote{This work was supported by Yandex Internet Company as well as by the grants of RFBR N 09-01-00294
and of the President of the Russian Federation N MD-8390.2010.1, NSH-691.2008.1.}}

\author{Liudmila Ostroumova
\footnote{Moscow State University, Mechanics and Mathematics Faculty, Department of Mathematical Statistics and Random Processes; Yandex, Departament of theoretical and practical research.}, 
Evgeniy Grechnikov
\footnote{Moscow State University, Mechanics and Mathematics Faculty, Department of Number Theory; Yandex, Departament of theoretical and practical research.}
}

\date{}

\renewcommand{\le}{\leqslant}
\renewcommand{\ge}{\geqslant}

\newcommand{\Prob}{{\rm P}}

\begin{document}

\newcommand{\MExpect}{\mathsf{M}}

\newenvironment{Proof}
{\par\noindent{\bf Proof}}
{\hfill$\scriptstyle$}

\maketitle

\section{Introduction}

In this paper we consider some properties of random graphs.
The standard random graph model $\mathfrak{G}(n,m)$ was introduced by Erd\H{o}s and R\'enyi in \cite{ER}. 
In this model we randomly choose one graph from all graphs with $n$ vertices and $m$ edges.
The similar model $\mathfrak{G}(n,p)$ was suggested by Gilbert in \cite{G}.
Here $n$ vertices are joined independently with probability $0<p<1$.
Many papers deal with the classical models. 
Main results can be found in \cite{B}, \cite{K}, \cite{JLR}.

Recently there has been interest in modeling complex real-world networks.
Real structures differ from standard random graphs.
One of the main characteristics of random graphs is their degree sequence.
In many real-world structures the degree sequence has a power law distribution. 
Standard random graph models do not have this property. 
So Barab\'asi and Albert suggested a new model in \cite{BA}. 
Then Bollob\'as and Riordan gave more precise definition of this model.
Many models of real-world networks and main results can be found in \cite{Res}.

This paper deals with the Bollob\'as--Riordan model.
Now let us describe this model.
Let $n$ be a number of vertices in our graph and $m$ be a fixed parameter. 
We begin with the case $m=1$.
We inductively construct a random graph $G_1^n$.
Start with $G_1^1$ the graph with one vertex and one loop.
Similarly we can start with $G_1^0$ the graph with no vertices.
Assume that we already constructed the graph $G_1^{t-1}$.
At the next step we add one vertex $t$ and one edge between vertices $t$ and $i$,
where $i$ is chosen randomly with
$$
\Prob(i=s) = \
\begin{cases}
d_{G_1^{t-1}}(s)/(2t-1) & \text{if } 1 \le s \le t-1,\cr \noalign{\smallskip}
{1/(2t-1)} &\text{if }  s=t. \cr
\end{cases}
$$ 

Here $d_{G_1^t}(s)$ is the degree of the vertex $s$ in $G_1^t$.
By $d(s)$ denote the degree of $s$ in the graph $G_1^n$.
In other words, the probability that a new vertex will be connected to the vertex $i$ is
proportional to the current degree of $i$.
Therefore this process is said to be {\it preferential attachment}. 
To obtain $G_m^n$ with $m>1$ we construct $G_1^{mn}$.
Then we identify the vertices $1, \dots, m$ to form the first vertex; we
identify the vertices $m+1, \dots, 2m$ to form the second vertex; and so on.
After this procedure, edges from $G_1^n$ connect ``big'' vertices in $G_m^n$.
Let $\mathfrak{G}_m^n$ be the probability space of constructed graphs.

Many papers deal with the Bollob\'as--Riordan model.
The diameter of this random graph was considered in \cite{BR1}.
In \cite{BR2} Bollob\'as and Riordan proved that
the degree sequence has a power law distribution.

\begin{Th}
If $m \ge 1$ is fixed, then there exists a function
$\varphi(n) = o(n)$ such that
for any $m \le d \le n^{1/15}$ we have 
$$
\lim_{n \to \infty} 
\Prob \left(\left|\#_m^n(d) - \frac{2nm(m+1)} {d(d+1)(d+2)}\right| >  \frac{\varphi(n)}{d(d+1)(d+2)}\right) = 0.
$$
Here $\#_m^n(d)$ is the number of vertices in $G_m^n$ with degree equal to $d$.
\end{Th}

Recently Grechnikov substantially improved Theorem 1 (see \cite{Gr}). 

In this paper we consider second degrees of vertices in $G_m^n$.
We estimate the expectation of the number of vertices with second degree equal to $d$.
Also we prove a concentration result.
This paper is organized as follows.
In section 2 we give main definitions and results.
In section 3 we prove all theorems.

\section{Definitions and results}

In this paper we study the random graph $ G_1^n $. 
When we write $ij \in G_1^n$ we mean that $G_1^n$ has the edge $ij$;
when we write $t \in G_1^n$ we simply mean that $t$ is a vertex of $G_1^n$.
Given a vertex $ t \in G_1^n $ we say that the {\it second degree} of the vertex $t$ is
$$
d_2(t) = \#\{ij: i \not= t, j \not= t, it \in G_1^n, ij \in G_1^n\}.
$$
In other words, the second degree of $t$ is the number of edges adjacent to the neighbors of $t$
except for the edges adjacent to the vertex $t$.

Let $M_n^1(d)$ be the expectation of the number of vertices with degree $d$ in $G_1^n$:
$$
M_n^1(d) = \MExpect\left(\#\{t \in G_1^n: d_{G_1^n}(t) = d\}\right).
$$
By  $X_n(d)$ denote the number of vertices with second degree $d$ in $G_1^n$.
By definition, put $M_n^2(d) = \MExpect X_n(d)$.

The aim of this paper is to prove the following results.

\begin{Th}
For any $k>1$ we have
$$
M_n^2(k) = \frac{4n}{k^2} \left(1 + O\left(\frac{\ln^2 k}{k}\right) + O\left(\frac{k^2}{n}\right) \right).
$$
\end{Th}

\begin{Th}
For any $\varepsilon > 0$ there exists a function $\varphi(n) = o(n)$ such that 
$$
\lim_{n \to \infty} \Prob \left(|X_n(k) - M_n^2(k)| \ge \frac{\varphi(n)}{k^2}\right) = 0
$$
for any $1 \le k \le n^{1/6 - \varepsilon}$.
\end{Th}
This is a concentration result which means that the distribution of second degrees does also obey 
(asymptotically) a power law.

To prove Theorem 2, we need the following definition.
Let $N_n(l,k)$ be the number of vertices in $G_1^n$ with degree $l$, with second degree $k$, and without loops:
$$
N_n(l,k) = \#\{t \in G_1^n: d(t)=l, d_2(t) = k, tt \notin G_1^n\}.
$$
We shall prove the following theorem.

\begin{Th}
In $G_1^n$ we have
$$
\MExpect N_n(l,k) = n \, c(l,k) \,(1+\theta(n,l,k)),
$$
where $|\theta(n,l,k)| < (2l+k-1)^2/n$.
The constants $c(l,k)$ are defined as follows:
\begin{eqnarray*}
c(l,0)&=& c(0,k)= 0,\\
c(1,k)&=&\frac{2k^2+14k}{(k+1)(k+2)(k+3)(k+4)},\\
c(l,k)&=&c(l,k-1)\frac{l+k-1}{2l+k+2}+c(l-1,k)\frac{l-1}{2l+k+2},\quad k>0,l>1.
\end{eqnarray*}
\end{Th}

We shall use the following lemmas to prove these theorems.

\begin{Lem}
Let $d \ge 1$ be natural; then
$$
M_n^1(d) = \frac {4 n} {d(d+1)(d+2)} \left(1 + \tilde \theta(n,d) \right),
$$
where $|\tilde \theta(n,d)| < d^2/n$. 
\end{Lem}

Denote by $P_n(l,k)$ the number of vertices in $G_1^n$ with a loop, with degree $l$, and with second degree $k$. 

\begin{Lem}
For any $n$ we have
$$
\MExpect P_n(l,k) \le p(l,k),
$$
where
\begin{eqnarray*}
p(2,0)&=&1,\\
p(l,k)&=&p(l,k-1)\frac{l+k-3}{2l+k-2}+p(l-1,k)\frac{l-1}{2l+k-2},\quad l\ge3,k\ge0.
\end{eqnarray*}
For the other values of $l$ and $k$ we have $p(l,k) = 0$.
\end{Lem}

The next section is organized as follows.
First we prove Theorem 4 and Theorem 2;
then we prove the lemmas.
Finally we give a proof of Theorem 3.

\section{Proofs}

\subsection{Proof of Theorem 4}

From the definition of $G_1^n$ it follows that $N_n(l,0) = N_n(0,k)=0$.
Indeed, since we have no vertices of degree $0$, we see that $N_n(0,k)=0$.
Since vertices with loops are not counted in $N_n(l,k)$,
it follows that we have no vertices of second degree $0$ and $N_n(l,0) =0$.
Therefore we have $\MExpect N_n(l,0)=\MExpect N_n(0,k)=0$.

Let us prove that
$\MExpect N_n(1,k) = n \, c(1,k) \,(1+\theta(n,1,k))$.
The proof is by induction on $k$.
For $k=0$ there is nothing to prove.
Now assume that for $j<k$ we have
$$
\MExpect N_n(1,j) = n \, c(1,j) \,(1+\theta(n,1,j)),
$$
where
$$|\theta(n,1,j)|<(j+1)^2/n,$$
$$c(1,j)=\frac {2j^2+14j} {(j+1)(j+2)(j+3)(j+4)}.$$

Denote by $N_i(l)$ the number of vertices with degree $l$ in $G_1^i$.

We need some additional notation.
Let $X$ be a function on $n$ (the number of vertices), $l$ (the first degree we are interested in), 
$k$ (the second degree we are interested in);
then denote by $\theta_1(X)$, $\theta_2(X)$, $\theta_3(X)$ $ \ldots $
some functions on $n$, $l$, $k$ such that
$|\theta_i(X)| < X$.

Obviously, $\MExpect N_{1}(1,k) = 0$.
For $i \ge 1$ we have
\begin{equation}\label{l=1}
\MExpect ( N_{i+1} (1,k) | N_{i}(1,k),N_{i}(1,k-1),N_{i}(k)) =  
N_i(1,k) \left( 1 - \frac {k+2} {2i+1} \right) + \frac {k N_i(1,k-1)} {2i+1} 
+ \frac {k\,N_i(k)} {2i+1}.
\end{equation}

Let us explain this equality.
Suppose we have $G_1^i$.
We add one vertex and one edge.
There are $N_i(1,k)$ vertices with degree $1$ and with second degree $k$ in $G_1^i$. 
The probability that we ``spoil'' one of these vertices is $(k+2) / (2i+1)$.
Also we have $N_i(1,k-1)$ vertices with degree $1$ and with second degree $k-1$. 
The probability that one of these vertices has degree $1$ and second degree $k$ in $G_1^{i+1}$ is $k/(2i+1)$. 
Finally, with probability equal to $k N_i(k)/(2i+1)$ the vertex $i+1$ has necessary degrees in $G_1^{i+1}$.

Using (\ref{l=1}), Lemma 1, and inductive assumption we get
$$
\MExpect N_{i+1} (1,k) =  
\MExpect N_i(1,k) \frac {2i-k-1} {2i+1}
+ \frac {k \MExpect N_i(1,k-1)} {2i+1}
+ \frac {k\,M_i^1(k)} {2i+1} =
$$
$$
=\MExpect N_i(1,k) \frac {2i-k-1} {2i+1}
+ \left(\frac {i\,k\,c(1,k-1)} {2i+1} 
+ \frac {4\,i} {(2i+1)(k+1)(k+2)}\right)
\left(1+\theta_1\left(k^2/i\right)\right).
$$
Let us introduce some notation:
$$
a_i = \frac {2i-k-1} {2i+1},
$$
$$
b_i =\frac {2i} {2i+1}
 \left(1+\theta_1\left(k^2/i\right)\right),
$$
$$
m = \frac {c(1,k-1)\,k} {2} + \frac {2} {(k+1)(k+2)}.
$$
Using this notation, we have 
$$
\MExpect N_{i+1} (1,k) =  
\MExpect N_i(1,k) \, a_i + m \, b_i.
$$

Let us prove the following equality by induction on $n$:
$$
\MExpect N_n(1,k) = \frac {2mn} {k+4} \left(1+\theta(n,1,k)\right).
$$
For $n=1$ we have $\MExpect N_{1}(1,k) = 0$. 
Since we have the condition $|\theta(1,1,k)|<(k+1)^2$,
we can take $\theta(1,1,k) = -1$.

Now put $t=k+1$. This is needed for the sequel.
Assume that
$$
\MExpect N_i(1,k) = \frac {2mi} {t+3} \left(1+\theta(i,1,t-1)\right).
$$
Then
$$
\MExpect N_{i+1} (1,k) =  
\MExpect N_i(1,k) \, a_i + m \, b_i =
$$
$$
=\frac {2mi(2i-t)} {(2i+1)(t+3)} \left(1+\theta_2(t^2/i)\right)
+ \frac {2mi} {2i+1}
 \left(1+\theta_1\left((t-1)^2/i\right)\right) =
$$
$$
= \frac {2m} {t+3}
\left( i+1 - \frac{1}{2i+1} + \theta_3\left(\frac{(2i-t)t^2}{2i+1} \right)
+\theta_4\left(\frac{(t-1)^2(t+3)}{2i+1} \right)\right).
$$
If $t \ge 1$ and $2i-t \ge 0$, then 
$$
\frac{1}{2i+1} + \frac{t^2|2i-t|}{2i+1}
+ \frac {(t-1)^2(t+3)}{2i+1} < t^2.
$$
Therefore,
$$
\MExpect N_{i+1} (1,k)
= \frac {2m(i+1)} {t+3}
\left( 1 + \theta_5\left(t^2/(i+1)\right)\right).
$$
In this case, we can put 
$\theta(i+1,1,k) = \theta_5\left(t^2/(i+1)\right)$.

If $t \ge 1$ and $2i - t \le -2$, then we do not have enough edges in $G_1^i$ and $\MExpect N_{i+1}(l,k) = 0$.
In this case, we can put $\theta(i+1,l,k) = -1$.

We consider the case $2i-t=-1$ later.

We get
$$
\MExpect N_n (1,k)
= \frac {2mn} {k+4}
\left( 1 + \theta(n,1,k)\right).
$$
Note that
$$
\frac {2m} {k+4} = \frac {4} {(k+1)(k+2)(k+4)} +  \frac {2\,c(1,k-1)\,k} {2(k+4)} =
$$
$$
= \frac {4} {(k+1)(k+2)(k+4)} +  \frac {2(k-1)^2+14(k-1)} {(k+1)(k+2)(k+3)(k+4)} =
$$
$$
= \frac {2k^2 +14k} {(k+1)(k+2)(k+3)(k+4)}
=c(1,k).
$$
This completes the proof for
$\MExpect N_n(1,k)$.

Consider the case $l,k>1$.
Assume that for all $i<l,j<k$ we have
$\MExpect N_n(i,j) = n \, c(i,j) \,(1+\theta(n,i,j))$.
Put $t = 2l+k-1$.
Obviously, $\MExpect N_{1}(l,k) = 0$.
For $i\ge 1$ we have
$$
\MExpect N_{i+1} (l,k) =  
\MExpect N_i(l,k) \left(1 - \frac {2l+k} {2i+1} \right) + \frac {(l-1)\,\MExpect N_i(l-1,k)} {2i+1}
+ \frac {(l+k-1) \MExpect N_i(l,k-1)} {2i+1} =
$$
$$
= \MExpect N_i(l,k) \frac {2i-t} {2i+1} +
 \left( \frac {(l-1)c(l-1,k)i} {2i+1} + \frac {(l+k-1)c(l,k-1)i} {2i+1} \right)
 \left( 1+\theta_1 \left((t-1)^2/i\right) \right). 
$$
Introduce some notation:
$$
a_i = \frac {2i-t} {2i+1}, 
$$
$$
b_i = \frac {2i} {2i+1} \left( 1+\theta_1\left((t-1)^2/i\right) \right),
$$
$$
m =\frac {(l-1)\,c(l-1,k)} {2} + \frac {(l+k-1)\,c(l,k-1)} {2}.
$$
We have
$$
\MExpect N_{i+1} (l,k) =  
\MExpect N_i(l,k) \, a_i + m \, b_i.
$$
It remains to prove the following statement by induction on $n$:
$$
\MExpect N_n(l,k) = 
\frac {2mn} {t+3} \left(1+\theta_5\left(t^2/n\right)\right) =
\frac {2mn} {t+3} \left(1+\theta(n,l,k)\right).
$$
The proof is the same as in the case of $l=1$.
In this case we have
$$
\frac {2m} {t+3} = 
\frac {(l-1)\,c(l-1,k)} {2l+k+2} + \frac {(l+k-1)\,c(l,k-1)} {2l+k+2}
=c(l,k).
$$

Now we need to consider only the case $2i-t=-1$. 
We need to show that $\MExpect N_{i+1}(l,k) = (i+1)c(l,k)(1+\theta(i+1,l,k))$.
We have $2(i+1) = 2l + k$.
In our graph $G_1^{i+1}$ we have $i+1$ edges.
Therefore the sum of all degrees is equal to $2l+k$.
Suppose we have at least one vertex with degree $l$ and second degree $k$.
We do not count vertices with a loop in $N_{i+1}(l,k)$.
Consequently $l$ edges go out from this vertex.
And there are $k/2$ edges between the neighbors of our vertex.
And we have no other edges.
Hence our vertex is joined to all other vertices in $G_1^{i+1}$.
So $l = i$. Thus $k=2$.
It follows that we consider the vertex 2.
And there is one edge from the vertex $2$ to the vertex $1$;
also edges from the vertices $3, \dots, i+1$ go to the vertex $2$.
So, there is only one graph with $N_{i+1}(l,k) \not= 0$.
This graph has only one vertex with degree $l$ and second degree $k$.
Therefore the probability of this graph is equal to $\MExpect N_{l+1}(l,2)$.
We have $\MExpect N_{l+1}(l,2) = \frac{2(l-1)!}{(2l+1)!!}$.

Recall that  $l=i$ and $k=2$.
Now we must only prove that 
$$
\MExpect N_{l+1} (l,2) =  
 (l+1)c(l,2)(1 + \theta(l+1,l,2)).
$$
Let us prove the inequality
$$
c(l,2) \ge \frac{24 (l-1)!}{5 (2l+4)!!}.
$$
It follows from the definition of $c(l,k)$ that
$$
c(1,2)=\frac {1} {10},
$$
$$
c(l,2) \ge c(l-1,2) \frac {l-1} {2l+4}, \,\,
l \ge 2.
$$

Obviously, $\theta(l+1,l,2) \ge -1$. 
Let us obtain the following upper bound:
$$
\theta(l+1,l,2) + 1 = \frac{\MExpect N_{l+1} (l,2)}{(l+1)c(l,k)} \le
\frac{2(l-1)! 5 (2l+4)!!}{(2l+1)!!(l+1) 24 (l-1)! } =
\frac{5 (2l+4)!!}{12 (2l+1)!!(l+1)} \le \frac{(2l+1)^2}{(l+1)}.
$$
This completes the proof. 

\subsection{Proof of Theorem 2}

From Theorem 4 we have the constants $c(l,k)$.
Imagine that we have a table with $c(l,k)$,
where $l$ is the number of a row and $k$ is the number of a column.
The sum of all numbers in the table is equal to $1$.
The sum of numbers in $l$-th row is equal to $\frac{4}{l(l+1)(l+2)}$.
It can easily be checked using the definition of  $c(l,k)$.
But we need to calculate $M_n^2(k)$,
so we are interested in the sum of all numbers in $k$-th column.
More precisely, 
$$
M_n^2(k) = \sum_{l=1}^{\infty} {\MExpect N_n(l,k)} + \sum_{l=1}^{\infty} {\MExpect P_n(l,k)}.
$$

First we estimate 
$\sum_{l=1}^{\infty} { c(l,k)}$.
Recall that
\begin{eqnarray*}
c(l,0)&=&0,\\
c(1,k)&=&\frac{2k^2+14k}{(k+1)(k+2)(k+3)(k+4)},\\
c(l,k)&=&c(l,k-1)\frac{l+k-1}{2l+k+2}+c(l-1,k)\frac{l-1}{2l+k+2},\quad k>0,l>1.
\end{eqnarray*}

Note that there exists a function $C(k) \ge 0$ such that for all $l\ge k\ge0$ and $l\ge1$ the inequality
\begin{equation}\label{cbound}
c(l,k)\le C(k)2^{-l}\frac{(l-1)!}{(l-k)!}
\end{equation}
holds. 
Indeed, the case of $k=0$ is obvious with $C(k)=0$.
In the case of $k\ge1$ we define $C(k)$ so that $C(k) \ge {C(k-1)}$ and (\ref{cbound}) holds for $l=k$.
We have
$$
(2l+k+2)\frac{c(l,k)}{C(k)}\le\frac{C(k-1)}{C(k)}(l+k-1)2^{-l}\frac{(l-1)!}{(l-k+1)!}+2^{-l+1}\frac{(l-1)!}{(l-k-1)!} \le
$$
$$
\le 2^{-l}\frac{(l-1)!}{(l-k)!}\left(\frac{l+k-1}{l-k+1}+2(l-k)\right)\le(2l+k+2)2^{-l}\frac{(l-1)!}{(l-k)!}.
$$
This proves (\ref{cbound}). 

In particular, the series $\sum_{l=1}^\infty l^N c(l,k)$ converges for
all $N$ and $k$.

Let us make some transformations:
$$
(2l+k+2)c(l,k)=(l+k-1)c(l,k-1)+(l-1)c(l-1,k),
$$
$$
\sum_{l=2}^\infty (2l+k+2)c(l,k)=\sum_{l=2}^\infty (l+k-1)c(l,k-1)+\sum_{l=1}^\infty lc(l,k),
$$
$$
\sum_{l=2}^\infty (l+k+2)c(l,k)=\sum_{l=2}^\infty (l+k-1)c(l,k-1)+c(1,k).
$$
Put $x_k=\sum_{l=2}^\infty c(l,k)$. 
Then $x_0=0$ and for $k\ge1$ we have
$$
(k+2)x_k=(k-1)x_{k-1}+c(1,k)+\sum_{l=2}^\infty l(c(l,k-1)-c(l,k)),
$$
$$
(k+2)(k+1)kx_k=(k-1)(k+1)kx_{k-1}+(k+1)kc(1,k)+\sum_{l=2}^\infty l(k(k+1)c(l,k-1)-k(k+1)c(l,k)),
$$
$$
(k+2)(k+1)kx_k=\sum_{s=1}^k(s(s+1)(s+2)x_s-(s-1)s(s+1)x_{s-1}) =
$$
$$=\sum_{s=1}^k s(s+1)c(1,s)+
\sum_{l=2}^\infty l\left(\sum_{s=1}^k(s(s+1)c(l,s-1)-s(s+1)c(l,s))\right)=
$$
$$
= \sum_{s=1}^k s(s+1)c(1,s)+\sum_{l=2}^\infty l\left(\sum_{s=1}^k((s+1)(s+2)-s(s+1))c(l,s)-(k+1)(k+2)c(l,k)\right)=
$$
$$
= \sum_{s=1}^k s(s+1)c(1,s)+\sum_{l=2}^\infty l\left(\sum_{s=1}^k 2(s+1)c(l,s)-(k+1)(k+2)c(l,k)\right),
$$
\begin{equation}\label{xbound}
x_k=\frac1{k(k+1)(k+2)}\sum_{s=1}^k s(s+1)c(1,s)+\frac2{k(k+1)(k+2)}\sum_{l=2}^\infty l\left(\sum_{s=1}^k(s+1)c(l,s)\right)
-\frac1k\sum_{l=2}^\infty lc(l,k).
\end{equation}

Put $y_k=\sum_{l=2}^\infty lc(l,k)$. 
Then
$$
x_k=\frac1{k(k+1)(k+2)}\sum_{s=1}^k s(s+1)c(1,s)+\frac2{k(k+1)(k+2)}\sum_{s=1}^k(s+1)y_s-\frac1k y_k.
$$
Make some transformations:
$$
(2l+k+2)lc(l,k)=(l+k-1)lc(l,k-1)+l(l-1)c(l-1,k),
$$
$$
\sum_{l=2}^\infty (2l+k+2)lc(l,k)=\sum_{l=2}^\infty(l+k-1)lc(l,k-1)+\sum_{l=1}^\infty l(l+1)c(l,k),
$$
$$
\sum_{l=2}^\infty(l+k+1)lc(l,k)=\sum_{l=2}^\infty(l+k-1)lc(l,k-1)+2c(1,k),
$$
$$
k y_k+\sum_{l=2}^\infty(l+1)lc(l,k)=(k-2)y_{k-1}+\sum_{l=2}^\infty l(l+1)c(l,k-1)+2c(1,k),
$$
$$
k(k-1)y_k=\sum_{s=1}^k(s(s-1)y_s-(s-1)(s-2)y_{s-1})=
$$
$$
=\sum_{s=1}^k\left((s-1)\sum_{l=2}^\infty l(l+1)c(l,s-1)-(s-1)\sum_{l=2}^\infty(l+1)lc(l,s)+2(s-1)c(1,s)\right)=
$$
$$
= 2\sum_{s=1}^k(s-1)c(1,s)+\sum_{l=2}^\infty l(l+1)\sum_{s=1}^kc(l,s)-k\sum_{l=2}^\infty l(l+1)c(l,k).
$$
For $k\ge2$ we have
$$
y_k=\frac2{k(k-1)}\sum_{s=1}^k(s-1)c(1,s)+\frac1{k(k-1)}\sum_{l=2}^\infty l(l+1)\sum_{s=1}^k c(l,s)-\frac1{k-1}\sum_{l=2}^\infty l(l+1)c(l,k).
$$
Let $z_k=\sum_{l=2}^\infty l(l+1)c(l,k)$. 
Then for $k\ge2$
$$
y_k=\frac2{k(k-1)}\sum_{s=1}^k(s-1)c(1,s)+\frac1{k(k-1)}\sum_{s=1}^k z_s-\frac1{k-1}z_k.
$$
Make similar transformations
$$
(2l+k+2)l(l+1)c(l,k)=(l+k-1)l(l+1)c(l,k-1)+(l+1)l(l-1)c(l-1,k),
$$
$$
\sum_{l=2}^\infty (2l+k+2)l(l+1)c(l,k)=\sum_{l=2}^\infty(l+k-1)l(l+1)c(l,k-1)+\sum_{l=1}^\infty l(l+1)(l+2)c(l,k),
$$
$$
\sum_{l=2}^\infty(l+k)l(l+1)c(l,k)=\sum_{l=2}^\infty(l+k-1)l(l+1)c(l,k-1)+6c(1,k),
$$
$$
\sum_{s=1}^k \sum_{l=2}^\infty(l+s)l(l+1)c(l,s)=\sum_{s=0}^{k-1} \sum_{l=2}^\infty(l+s)l(l+1)c(l,s)+\sum_{s=1}^k 6c(1,s),
$$
$$
\sum_{l=2}^\infty(l+k)l(l+1)c(l,k)=\sum_{s=1}^k 6c(1,s).
$$
Since $c(1,s)=O\left(\frac1{s^2}\right)$, we have
$$
0\le z_k\le\frac1k\sum_{l=2}^\infty(l+k)l(l+1)c(l,k)=O\left(\frac1k\sum_{s=1}^k\frac1{s^2}\right)=O\left(\frac1k\right),
$$
$$
\sum_{s=1}^k(s-1)c(1,s)=O\left(\sum_{s=1}^k\frac1s\right)=O(\ln k),
$$
$$
y_k=O\left(\frac{\ln k}{k^2}\right),
$$
$$
\sum_{s=1}^k(s+1)y_s=O\left(\sum_{s=1}^k\frac{\ln s}{s}\right)=O(\ln^2k),
$$
$$
x_k=\frac1{k(k+1)(k+2)}\sum_{s=1}^k s(s+1)c(1,s)+O\left(\frac{\ln^2k}{k^3}\right).
$$
Finally, $c(1,s)=\frac2{s(s+1)}+O\left(\frac1{s^3}\right)$, so $\sum_{s=1}^k s(s+1)c(1,s)=2k+O(\ln k)$ and
$$
x_k=\frac2{(k+1)(k+2)}+O\left(\frac{\ln^2k}{k^3}\right)=\frac2{k^2}+O\left(\frac{\ln^2k}{k^3}\right),
$$
$$
\sum_{l=1}^\infty c(l,k)=c(1,k)+x_k=\frac4{k^2}+O\left(\frac{\ln^2k}{k^3}\right).
$$

Now we can estimate $M_n^2(k)$:
$$
M_n^2(k) = \sum_{l=1}^{\infty} {c(l,k)\,n(1 + \theta(n,l,k))} + \sum_{l=1}^{\infty} {\MExpect P_n(l,k)}.
$$
The first sum:
$$
 \sum_{l=1}^{\infty} {c(l,k)\,n} 
 = \frac{4n}{ k^2}+O\left(\frac{n \ln^2k}{k^3}\right).
$$
The second sum:
$$
\sum_{l=1}^{\infty}{c(l,k)\,n\,|\theta(n,l,k)|} \le
\sum_{l=1}^{\infty}{c(l,k)\,(2l+k)^2} =
\sum_{l=1}^{\infty} 4 l^2 c(l,k)  + \sum_{l=1}^{\infty} 4lk c(l,k) + \sum_{l=1}^{\infty} k^2 c(l,k)  =
$$
$$
= 4 c(1,k) + \sum_{l=2}^{\infty} 4 l(l+1) c(l,k) - \sum_{l=2}^{\infty} 4 l c(l,k) + 
4 k c(1,k) + \sum_{l=2}^{\infty} 4lk c(l,k) +
\sum_{l=1}^{\infty} k^2 c(l,k) =
$$
$$
= (4 + 4k) c(1,k) + 4 z_k + (4 k - 4) y_k + \sum_{l=1}^{\infty} k^2 c(l,k)
= O\left(\frac{1}{k} + \frac{1}{k} + \frac{\ln k}{k} + \frac{\ln^2 k}{k} + 1 \right)
=O\left(1\right).
$$
The third sum:
$$
\sum_{l=1}^{\infty} {\MExpect P_n(l,k)} \le
\sum_{l=1}^{\infty} p(l,k).
$$
Recall that
\begin{eqnarray*}
p(2,0)&=&1,\\
p(l,k)&=&p(l,k-1)\frac{l+k-3}{2l+k-2}+p(l-1,k)\frac{l-1}{2l+k-2},\quad k\ge0,l\ge3.
\end{eqnarray*}
For the other values of $l$ and $k$ we have $p(l,k) = 0$.
We can estimate $p(l,k)$:
$$
p(l,k) \le \frac{6}{l(l+1)}. 
$$
Indeed, it is easy to check that the function $\frac{6}{l(l+1)}$ follows the recurrent relation.
So when $l=2$ and $k=0$ we use the fact that $p(l,k) = 1 \le \frac{6}{l(l+1)}$, and then we proceed by induction.
Hence the series $\sum_{l=2}^\infty p(l,k)$ converges. In other words, $\sum_{l=2}^\infty p(l,k) = O(1)$. 
Therefore
$$
M_n^2(k) = \frac{4n}{k^2}+O\left(\frac{n \ln^2k}{k^3}\right)+ O(1)+O(1) =
\frac{4n}{k^2}\left(1 + O\left(\frac{\ln^2k}{k}\right)+ O\left(\frac{k^2}{n}\right)\right).
$$
This completes the proof. 

Now we must only prove Lemma 1 and Lemma 2.

\subsection{Proof of Lemma 1}
In \cite{BR2} Bollob\'as and Riordan computed the expectation of the number of vertices with degree $d$.
But they only looked at $d \le n^{1/15}$ and proved that
$$
\MExpect M_n^1(d) \sim \frac{4 n}{d(d+1)(d+2)}.
$$
We are interested in $\MExpect M_n^1(d)$ for any $d$.
In addition, we want to estimate $|\tilde\theta(n,d)|$.
Therefore we compute $M_n^1(d)$ in this paper.

The proof is by induction on $d$.
First we need to consider 2 cases: $d = 1$ and $d = 2$.

Consider the case $d=1$.
Obviously, $M_0^1(1) = 0$.
Assume that
$M_i^1(1) = \frac {2i} {3} \left(1 + \tilde\theta(i,1)\right)$.
Then
$$
M_{i+1}^1 (1) =  
M_i^1(1)\left(1 - \frac{1}{2i+1}\right) + \frac{2i}{2i+1} = 
\frac {2i} {3} \left(1 + \tilde\theta(i,1)\right)\frac{2i}{2i+1} + \frac{2i}{2i+1} =
$$
$$
=\frac{2}{3} \left(i+1 - \frac{1}{2i+1} +\frac{2i^2}{2i+1} \tilde\theta(i,1)\right) =
\frac{2(i+1)}{3} \left(1 - \frac{1}{(2i+1)(i+1)} + \frac{2i^2}{(2i+1)(i+1)} \tilde\theta(i,1)\right).
$$
Put 
$\tilde\theta(i+1,1) = \frac{2i^2}{(2i+1)(i+1)} \tilde\theta(i,1)- \frac{1}{(2i+1)(i+1)}$.
Note that
$$
|\tilde\theta(i+1,1)| \le 
\frac{2i}{(2i+1)(i+1)} +  \frac{1}{(2i+1)(i+1)} \le
1/(i+1).
$$
This completes the proof for $d=1$.

The case $d=2$ is somewhat different.
Obviously, $M_0^1(2) = 0$.
Suppose
$M_i^1(2) = \frac {i} {6} \left(1 + \tilde\theta(i,2)\right)$.
Then
$$
M_{i+1}^1 (2) =  
M_i^1(2)\left(1 - \frac{2}{2i+1}\right) + M_i^1(1)\frac{1}{2i+1} + \frac{1}{2i+1} = 
$$
$$
=\frac {i} {6} \left(1 + \tilde\theta(i,2)\right)\frac{2i-1}{2i+1} + 
\frac {2i} {3(2i+1)} \left(1 + \tilde\theta(i,1)\right) + \frac{1}{2i+1} =
$$
$$
=\frac{1}{6} \left(i+1 + \frac{5}{2i+1} +
\frac{(2i-1)i}{2i+1}\tilde\theta(i,2) + \frac{4i}{2i+1}\tilde\theta(i,1)\right) =
$$
$$
=\frac{i+1}{6} \left(1 + \frac{5}{(2i+1)(i+1)} + 
\frac{(2i-1)i}{(2i+1)(i+1)}\tilde\theta(i,2) + \frac{4i}{(2i+1)(i+1)}\tilde\theta(i,1)\right).
$$
Put
$$
\tilde\theta(i+1,2) = 
\frac{5}{(2i+1)(i+1)} + \frac{(2i-1)i}{(2i+1)(i+1)}\tilde\theta(i,2) + \frac{4i}{(2i+1)(i+1)}\tilde\theta(i,1).
$$
Note that $\tilde\theta(i,1)<0$.
Hence
$$
|\tilde\theta(i+1,2)| \le \left|\frac{(2i-1)i}{(2i+1)(i+1)}\tilde\theta(i,2)\right| + 
\max \left\{ \left|\frac{5}{(2i+1)(i+1)}\right|,
\left|\frac{4i}{(2i+1)(i+1)}\tilde\theta(i,1)\right|
\right\}.
$$
We got necessary bounds for
$\tilde\theta(i,2)$ and $\tilde\theta(i,1)$. 
Thus, it is easy to check that 
$$
|\tilde\theta(i+1,2)| \le \frac {4} {i+1}.
$$
This completes the proof for $d=2$.

Suppose $d \ge 3$ and we can prove the theorem for all smaller degrees.
This case is proved by induction on $i$.
For $i=0$ we have $M_0^1(d) = 0$.
Assume that
$M_i^1(d) = \frac {4i} {d(d+1)(d+2)} \left(1 + \tilde\theta(i,d)\right)$.
Then
$$
M_{i+1}^1 (d) =  
M_i^1(d)\left(1 - \frac{d}{2i+1}\right) + M_i^1(d-1)\frac{d-1}{2i+1} = 
$$
$$
= \frac {4i(2i+1-d)} {d(d+1)(d+2)(2i+1)}\left(1 + \tilde\theta(i,d)\right)
+ \frac{4i}{d(d+1)(2i+1)}\left(1 + \tilde\theta(i,d-1)\right) =
$$
$$
= \frac{4(i+1)}{d(d+1)(d+2)} \left(1 - \frac{1}{(2i+1)(i+1)} +
\frac{i(2i+1-d)}{(2i+1)(i+1)}\tilde\theta(i,d) + \frac{i(d+2)}{(2i+1)(i+1)}\tilde\theta(i,d-1) \right).
$$
If $2i+1-d \ge 0$, we can put
$$
\tilde\theta(i+1,d) = 
- \frac{1}{(2i+1)(i+1)} +
\frac{i(2i+1-d)}{(2i+1)(i+1)}\tilde\theta(i,d) + \frac{i(d+2)}{(2i+1)(i+1)}\tilde\theta(i,d-1).
$$
We obtain the following estimate:
$$
|\tilde\theta(i+1,d)| \le
\frac{1}{(2i+1)(i+1)} +
\frac{i(2i+1-d)}{(2i+1)(i+1)}|\tilde\theta(i,d)| + \frac{i(d+2)}{(2i+1)(i+1)}|\tilde\theta(i,d-1)| \le
\frac{d^2}{i+1}.
$$
If $2i+2-d \le 0$, we have no vertices with degree $d$ in $G_1^{i+1}$. 
Indeed, in $G_1^{i+1}$ the sum of all degrees is $2i+2$.
If $d > 2i + 2$, we obviously do not have enough edges.
If $d = 2i + 2$, then it is easy to check that we can not have any vertices with degree $d$ ($d > 2$). 
So we can put $\tilde\theta(i+1,d) = -1$.
This concludes the proof.

\subsection{Proof of Lemma 2}

Obviously,
$\MExpect P_n(0,k)= \MExpect P_n(1,k)=0$. 
For all $k>0$ we have $\MExpect P_n(2,k)=0$.
For $k=0$ we have
$$
\MExpect P_n(2,0) = \sum_{i=1}^n {\frac {1} {2i-1} \prod_{j=i+1}^{n} \frac {2j-3} {2j-1}}
= \sum_{i=1}^n {\frac{1}{2n-1}}
= \frac {n} {2n-1} 
\le 1.
$$

The rest of the proof is by induction.
Consider $l\ge3$, $k\ge0$.
Suppose that for $i<l$ and $j<k$ we have
$\MExpect P_n(i,j) \le p(i,j)$.

Trivially, $P_1(l,k) = 0$.
It is easily shown that 
$\MExpect P_{i+1}(l,k) = 0$ if $2i + 4 < 2l+k$.

If $2i + 4 = 2l + k$ and $P_{i+1}(l,k) \not= 0$, then $l = i+2$ and $k = 0$.
And we have only one graph with $P_{i+1}(l,k) \not= 0$.
Arguing as in the end of Section 3.1, we see that the probability of this graph is 
$\frac {(l-1)!}{(2l-1)!!}$.
From the recurrent relation we have $p(l,0) = \frac{1}{2^{l-2}}$.
In our case we get
$$
\MExpect P_{i+1}(l,k) = \frac {(l-1)!}{(2l-1)!!} < \frac{1}{2^{l-2}} = p(l,0).
$$

If $2i + 3 \ge 2l + k$, then
$$
\MExpect P_{i+1} (l,k) =  
\MExpect P_i(l,k)\left(1-\frac{2l+k-2}{2i+1}\right) + \MExpect P_i(l,k-1)\frac{l+k-3}{2i+1}
+\MExpect P_i(l-1,k)\frac{l-1}{2i+1}.
$$
Using the recurrent relation for $p(l,k)$ and induction on $i$ it is easy to prove that 
$\MExpect P_n(l,k) \le p(l,k)$.
This concludes the proof of Lemma 2.

\subsection{Proof of Theorem 3}

This proof is similar to the proof given in \cite{BR2}.
But our case is more complicated.
We need 
the Azuma--Hoeffding inequality (see \cite{A}):

\begin{Lem}
Let $(X_i)_{i=0}^{n}$ be a martingale such that
$|X_i - X_{i-1}| \le c$ for any $1 \le i \le n$.
Then
$$
\Prob\left(|X_n - X_0| \ge x \right) \le 2e^{-\frac{x^2}{2c^2n}}
$$
for any $x>0$.
\end{Lem}

Suppose we are given an $\varepsilon > 0$.
Fix $n \ge 3$ and $k$: $1 \le k \le n^{1/6-\varepsilon}$.
Consider the random variables
$X^i(k) = \MExpect(X_n(k)|G_1^i)$,
$i = 0, \ldots, n$.
Let us explain the notation $\MExpect(X_n(k)|G_1^i)$.
We construct the graph $G_1^n \in \mathfrak{G}_1^n$ by induction.
For any $t \le n$ there exists a unique
$G_1^t \in \mathfrak{G}_1^t$ such that 
$G_1^n$ is obtained from $G_1^t$.
So $\MExpect(X_n(k)|G_1^t)$ is the expectation of the number of vertices with second degree $k$ in $G_1^n$
if at the step $t$ we have the graph $G_1^t$.

Note that
$X^0(k) = \MExpect X_n(k)$ and $X^n(k) = X_n(k)$.
From the definition of $G_1^n$ it follows that $X^i(k)$ is a martingale.

We will prove below that for any $i = 1, \ldots, n$
$$
|X^{i}(k) - X^{i-1}(k)| \le 10\,k\,\ln{n}.
$$
Theorem 3 follows from this statement immediately.
Put $c = 10\,k\ln{n}$.
Then from Azuma--Hoeffding inequality it follows that
$$
\Prob\left(|X_n(k) - \MExpect X_n(k)| \ge k\, \sqrt{n} \, \ln^2{n}\right) \le
2\exp\left\{-\frac{n \,k^2 \,\ln^4{n} }{200\,n\, k^2\ln^2{n}} \right\} = o(1).
$$
If $k \le n^{1/6-\varepsilon}$, then the value of
$n/k^2$ is considerably greater than $k \, \ln^2{n}\, \sqrt{n}$.
This means that we have
$\left(k \, \sqrt{n}\, \ln^2{n}\right)/\left(n/k^2  \right) = o(1)$.
This is exactly what we need.

It remains to estimate the quantity $|X^{i}(k) - X^{i-1}(k)|$.
The proof is by a direct calculation.

Fix $1 \le i \le n$ and some graph $G_1^{i-1}$.
Note that
$$
\left|\MExpect\left(X_n(k)|G_1^{i}\right) - \MExpect\left(X_n(k)|G_1^{i-1}\right)\right| \le
\max_{\tilde G_1^i\supset G_1^{i-1}}  \left\{\MExpect\left(X_n(k)|\tilde G_1^{i}\right)\right\} -
\min_{\tilde G_1^i\supset G_1^{i-1}}  \left\{\MExpect\left(X_n(k)|\tilde G_1^{i}\right)\right\}.
$$

Put
$\hat G_1^i = \arg \max \MExpect(X_n(k)|\tilde G_1^{i})$,
$\bar G_1^i = \arg \min \MExpect(X_n(k)|\tilde G_1^{i})$.
We need to estimate the difference
$\MExpect(X_n(k)|\hat G_1^{i}) - \MExpect(X_n(k)|\bar G_1^{i})$.

Using the notation $N_n(l,k)$ and $P_n(l,k)$ from Section 2, we get
$$
\MExpect(X_n(k)|\hat G_1^{i}) =
\sum_{l=1}^{\infty} \MExpect(N_n(l,k)|\hat G_1^{i})
+ \sum_{l=1}^{\infty} \MExpect(P_n(l,k)|\hat G_1^{i}),
$$
$$
\MExpect(X_n(k)|\bar G_1^{i}) =
\sum_{l=1}^{\infty} \MExpect(N_n(l,k)|\bar G_1^{i})
+ \sum_{l=1}^{\infty} \MExpect(P_n(l,k)|\bar G_1^{i}).
$$

For $i \le t \le n$ put 
$$
\delta_t(l,k) = \MExpect(N_t(l,k)|\hat G_1^{i}) - \MExpect(N_t(l,k)|\bar G_1^{i}),
\,\,\,\,
\delta'_t(l,k) = \delta_t(l,k) I(\delta_t(l,k) > 0),
$$
$$
\epsilon_t(l,k) = \MExpect(P_t(l,k)|\hat G_1^{i}) - \MExpect(P_t(l,k)|\bar G_1^{i}),
\,\,\,\,
\epsilon'_t(l,k) = \epsilon_t(l,k) I(\epsilon_t(l,k) > 0),
$$
$$
\delta_t(k) = \MExpect{(N_t(k)|\hat G_1^{i})} - 
\MExpect(N_t(k)|\hat G_1^{i})),
\,\,\,\,
\delta'_t(k) = \delta_t(k) I(\delta_t(k) > 0).
$$

Note that
$$
\MExpect(X_n(k)|\hat G_1^{i}) - \MExpect(X_n(k)|\bar G_1^{i}) =
\sum_{l=1}^{\infty} \delta_n(l,k)
+ \sum_{l=1}^{\infty} \epsilon_n(l,k) \le
$$
$$
\le \sum_{l=1}^{\infty} \delta'_n(l,k)
+ \sum_{l=1}^{\infty} \epsilon'_n(l,k) \le
\sum_{l=1}^{\infty} \sum_{j=0}^{k} \left( \delta'_n(l,j) + \epsilon'_n(l,j) \right).
$$
Let us estimate this double sum.

First suppose that $n=i$. 
Fix $G_1^{i-1}$.
Graphs $\hat G_1^i$ and $\bar G_1^i$ are obtained from the graph $G_1^{i-1}$.
We add the vertex $i$ and one edge $i\hat q$ or $i\bar q$, respectively.
New edge changes only the degree of $\hat q$ or $\bar q$
and the second degree of neighbors of $\hat q$ or $\bar q$, respectively.
Consider $\hat G_1^i$.
Fix $l$ and $j \le k$.
We are interested in measuring the growth of the number of vertices with  degree $l$ and second degree $j$ at the step $i$. 
First $i$ can become a vertex of second degree $j$ with $j \le k$.
Secondly the vertex $\hat q$ can become a vertex of second degree $j$ with $j \le k$.
Thirdly  the second degree of neighbors of $\hat q$ increases.
If $\hat q$ has at least $k+1$ neighbors in $G_1^{i-1}$, 
then after the step $i$ these vertices have second degree bigger than $k$ and we do not count them. 
If $\hat q$ has at most $k$ neighbors in $G_1^{i-1}$, 
then at most $k$ vertices change their second degrees at the step $i$. 
Arguing as above, we consider $\bar G_1^{i}$.
We are interested in measuring the decrease of the values $N_{i-1}(l,j)$ and $P_{i-1}(l,j)$.
First $\bar q$ has new degree after the step $i$. 
Secondly some neighbors of $\bar q$ can have second degree $j \le k$ in $G_1^{i-1}$ 
(so the number of the neighbors of $\bar q$ in $G_1^{i-1}$ is not bigger than $k+1$).
Let us sum all the just-mentioned numbers.
We have
$$
\sum_{l=1}^{\infty} \sum_{j=0}^{k} \left( \delta'_i(l,j) + \epsilon'_i(l,j) \right) \le 
1+1+k+1+(k+1) = 2k+4.
$$

The case $n=i$ is complete.
Now consider $t$: $i \le t \le n-1$.
Note that
$$
\MExpect \left(N_{t+1} (1) |G_1^i\right) =  
\MExpect\left(N_{t} (1) |G_1^i\right)\left(1 - \frac{1}{2t+1}\right) + \frac{2t}{2t+1},
$$
$$
\MExpect\left(N_{t+1} (2) |G_1^i\right) =  
\MExpect\left(N_{t} (2) |G_1^i\right)\left(1 - \frac{2}{2t+1}\right) + 
\MExpect\left(N_{t} (1) |G_1^i\right)\frac{1}{2t+1} + \frac{1}{2t+1},
$$
$$
\MExpect\left(N_{t+1} (j) |G_1^i\right) =  
\MExpect\left(N_{t} (j) |G_1^i\right)\left(1 - \frac{j}{2t+1}\right) + 
\MExpect\left(N_{t} (j-1) |G_1^i\right)\frac{j-1}{2t+1}, \,\,\, j \ge 3,
$$
$$
\MExpect \left(N_{t+1} (1,j) |G_1^i\right) =  
\MExpect\left(N_t(1,j)|G_1^i\right) \left( 1 - \frac {j+2} {2t+1} \right) + 
\frac {j \, \MExpect\left(N_t(1,j-1)|G_1^i\right)} {2t+1} 
+ \frac {j\,\MExpect\left(N_t(j)|G_1^i\right)} {2t+1},
$$
$$
\MExpect \left(N_{t+1} (l,j) |G_1^i\right)=  
\MExpect \left( N_t(l,j)|G_1^i\right) \left(1 - \frac {2l+j} {2t+1} \right) + 
\frac {(l-1)\,\MExpect \left( N_t(l-1,j)|G_1^i\right)} {2t+1} +
$$
$$
+ \frac {(l+j-1) \MExpect \left(N_t(l,j-1)|G_1^i\right)} {2t+1}, \,\,\, l \ge 2,
$$
$$
\MExpect \left(P_{t+1} (2,0) |G_1^i\right) =  
\MExpect \left(P_{t} (2,0) |G_1^i\right)\left(1-\frac{2}{2t+1}\right)
+\frac{1}{2t+1},
$$
$$
\MExpect \left(P_{t+1} (l,j) |G_1^i\right) =  
\MExpect \left(P_{t} (l,j) |G_1^i\right)\left(1-\frac{2l+j-2}{2t+1}\right) + 
\MExpect \left(P_{t} (l,j-1) |G_1^i\right)\frac{l+j-3}{2t+1}+
$$
$$
+\MExpect \left(P_{t} (l-1,j) |G_1^i\right)\frac{l-1}{2t+1}, \,\,\, l \ge 3.
$$ 
We obtained the same equalities in proofs of Theorem 4, Lemma 1, and Lemma 2.
Replace $G_1^i$ by $\hat G_1^i$ or $\bar G_1^i$ in these equalities.
Substracting the equalities with $\bar G_1^i$  from the equalities with $\hat G_1^i$
and using the inequality $(a+b)I(a+b>0) \le aI(a>0) + bI(b>0)$, we get
$$
\delta'_{t+1}(j) \le \delta'_t(j) \left(1 - \frac{j}{2t+1} \right) + \delta'_t(j-1)\frac{j-1}{2t+1},
$$
$$
\delta'_{t+1}(1,j) \le \delta'_t(1,j) \left(1 - \frac{j+2}{2t+1}\right) + \frac{j\delta'_t(1,j-1)}{2t+1}
+ \frac{j\delta'_t(j)}{2t+1},
$$
$$
\delta'_{t+1}(l,j) \le \delta'_t(l,j) \left(1 - \frac{2l+j}{2t+1}\right) + \frac{(l-1)\delta'_t(l-1,j)}{2t+1} +
\frac{(l+j-1)\delta'_t(l,j-1)}{2t+1}, \,\,\, l \ge 2,
$$
$$
\epsilon'_{t+1}(2,0) \le \epsilon'_{t}(2,0) \left( 1 - \frac{2}{2t+1} \right),
$$
$$
\epsilon'_{t+1}(l,j) \le \epsilon'_t(l,j) \left(1 - \frac{2l+j-2}{2t+1}\right) + \frac{(l-1)\epsilon'_t(l-1,j)}{2t+1} +
\frac{(l+j-3)\epsilon'_t(l,j-1)}{2t+1}, \,\,\, l \ge 3.
$$
Now we can estimate the sum
$$
\sum_{l=1}^{\infty} \sum_{j=0}^{k} \delta'_{t+1}(l,j)
+ \sum_{l=2}^{\infty} \sum_{j=0}^{k} \epsilon'_{t+1}(l,j) \le
$$
$$
 \le \sum_{j=1}^{k} \left( \delta'_t(1,j) \left(1 - \frac{j+2}{2t+1}\right) + \frac{j\delta'_t(1,j-1)}{2t+1}
 + \frac{j\delta'_t(j)}{2t+1} \right) +
 \epsilon'_{t}(2,0) \left( 1 - \frac{2}{2t+1} \right)+
$$
$$
+ \sum_{l=2}^{\infty} \sum_{j=1}^{k} \left(\delta'_t(l,j) \left(1 - \frac{2l+j}{2t+1}\right) + \frac{(l-1)\delta'_t(l-1,j)}{2t+1} +
\frac{(l+j-1)\delta'_t(l,j-1)}{2t+1} \right) +
$$
$$
+ \sum_{l=3}^{\infty} \sum_{j=0}^{k} \left( \epsilon'_t(l,j) \left(1 - \frac{2l+j-2}{2t+1}\right) + \frac{(l-1)\epsilon'_t(l-1,j)}{2t+1} +
\frac{(l+j-3)\epsilon'_t(l,j-1)}{2t+1} \right) =
$$

$$
=
\sum_{j=1}^{k} \delta'_t(1,j) 
- \sum_{j=1}^{k} \frac{(j+2)\delta'_t(1,j)}{2t+1} 
+ \sum_{j=1}^{k-1} \frac{(j+1)\delta'_t(1,j)}{2t+1}
+ \sum_{j=1}^{k} \frac{j\delta'_t(j)}{2t+1} + \sum_{l=2}^{\infty} \sum_{j=1}^{k} \delta'_t(l,j) -
$$
$$  
- \sum_{l=2}^{\infty} \sum_{j=1}^{k} \frac{(2l+j)\delta'_t(l,j)}{2t+1}
+ \sum_{l=1}^{\infty} \sum_{j=1}^{k} \frac{l\delta'_t(l,j)}{2t+1}
+ \sum_{l=2}^{\infty} \sum_{j=1}^{k-1} \frac{(l+j)\delta'_t(l,j)}{2t+1} +
$$
$$
+ \epsilon'_{t}(2,0)  - \frac{2 \epsilon'_{t}(2,0)}{2t+1} 
+ \sum_{l=3}^{\infty} \sum_{j=0}^{k} \epsilon'_t(l,j)
- \sum_{l=3}^{\infty} \sum_{j=0}^{k} \frac{(2l+j-2)\epsilon'_t(l,j) }{2t+1} + 
$$
$$
+ \sum_{l=3}^{\infty} \sum_{j=0}^{k} \frac{l\epsilon'_t(l,j)}{2t+1}
+ \frac{2\epsilon'_t(2,0)}{2t+1}
+ \sum_{l=3}^{\infty} \sum_{j=0}^{k-1} \frac{(l+j-2)\epsilon'_t(l,j)}{2t+1} =
$$

$$
= \sum_{l=1}^{\infty} \sum_{j=0}^{k} \left( \delta'_{t}(l,j)
+ \epsilon'_{t}(l,j) \right)
+ \sum_{j=1}^{k} \frac{j\delta'_t(j)}{2t+1}
- \sum_{l=1}^{\infty} \delta'_t(l,k)\frac{l+k}{2t+1}
- \sum_{l=3}^{\infty} \epsilon'_t(l,k)\frac{l+k-2}{2t+1} \le
$$
$$
\le \sum_{l=1}^{\infty} \sum_{j=0}^{k} \left( \delta'_{t}(l,j)
+ \epsilon'_{t}(l,j) \right) + \sum_{j=1}^{k} \frac{j\delta'_t(j)}{2t+1}.
$$
It remains to estimate the sum $\sum_{j=1}^{k} \frac{j\delta'_t(j)}{2t+1}$.
Note that for any $t \ge i$ we have $\sum_{j=0}^{k} \delta'_t(j) \le 3$.
It is obvious for $t=i$ 
(when we add a new vertex $i$, we change only the degree of $\hat q$ or $\bar q$).
If $t+1>i$, then
$$
\sum_{j=1}^{k} \delta'_{t+1}(j) \le 
\sum_{j=1}^{k} \left( \delta'_t(j) \left(1 - \frac{j}{2t+1} \right) + \delta'_t(j-1)\frac{j-1}{2t+1} \right)
=\sum_{j=1}^{k} \delta'_{t}(j) - \delta'_t(k)\frac{k}{2t+1}.
$$
In other words,  $\sum_{j=1}^{k} \delta'_t(j)$ is not increasing when $t$ is growing.

So we get
$$
\sum_{l=1}^{\infty} \sum_{j=0}^{k} \left( \delta'_{t+1}(l,j)
+ \epsilon'_{t+1}(l,j) \right) \le
\sum_{l=1}^{\infty} \sum_{j=0}^{k} \left( \delta'_{t}(l,j)
+ \epsilon'_{t}(l,j) \right) + \frac{3 k}{2t+1}.
$$
Thus we have
$$
|\MExpect(X_n(k)|G_1^{i-1}) - \MExpect(X_n(k)|G_1^{i})| \le
\sum_{l=1}^{\infty} \sum_{j=0}^{k} \left( \delta'_{i}(l,j)
+ \epsilon'_{i}(l,j) \right) + \sum_{t=i}^{n-1}\frac{3 k}{2t+1} \le
$$
$$
\le 2k+4+ \sum_{t=1}^{n-1}\frac{3 k}{2t+1} \le
2k+5 + \frac{3}{2} k\ln n \le 10\,k \ln n.
$$
This concludes the proof of Theorem 3.

\end{document}